\documentclass[12pt]{article}
\usepackage{amsmath,amssymb,amscd,latexsym,euscript,mathrsfs}
\usepackage[matrix,arrow,curve]{xy}

\newcommand{\dom}{ {\rm dom\,}}
\newcommand{\cod}{ {\rm cod\,}}
\newcommand{\coLim}{\underrightarrow{\lim}}

\newcommand{\Set}{{\rm Set}}
\newcommand{\Ab}{{\rm Ab}}
\newcommand{\Ob}{{\rm Ob\,}}
\newcommand{\Mor}{{\rm Mor\,}}

\newcommand{\mB}{{\mathscr B}}
\newcommand{\mC}{{\mathscr C}}
\newcommand{\mD}{{\mathscr D}}

\newcommand{\ZZ}{{\,\mathbb Z}}

\newcommand{\mA}{{\mathcal A}}

\newcommand{\fF}{{\mathfrak{F}}}

\newtheorem{theorem}{\bf Theorem}[section]
\newtheorem{lemma}[theorem]{\bf Lemma}
\newtheorem{proposition}[theorem]{\bf Proposition}
\newtheorem{corollary}[theorem]{\bf Corollary}
\newtheorem{definition}{\sc Definition}[section]

\def\leq{\leqslant}
\def\geq{\geqslant}

\begin{document}

\begin{center}
 {\large THE HOMOLOGY GROUPS OF A PARTIAL TRACE MONOID ACTION
}
\end{center}
\centerline{A. A. Husainov, husainov51@yandex.ru}

\begin{abstract}
The aim of this paper is to investigate the homology groups of mathematical models of 
concurrency. 
We study the Baues-Wirsching homology groups of a small category
 associated with a partial monoid action on a set.
We prove that these groups can be reduced to the Leech homology groups of the monoid.
For a trace monoid with an action on a set, we will build
a cubical complex of free Abelian groups with homology groups isomorphic to the integral homology groups
of the action category.
It allows us to solve the
  problem posed by the author in 2004 of the constructing an algorithm for
computing homology groups of the CE nets. 
We describe the algorithm and give examples of calculating the homology groups.
\end{abstract}

2000 Mathematics Subject Classification 18B40, 18G10, 18G35, 55U10, 68Q85


\section*{Introduction}

The homology groups of a partial trace monoid 
action on a set was studied in \cite{X20042}.
They were used to introduce the homology groups of CE nets.
 An algorithm for computing the first homology group of CE nets was obtained.
The question of constructing an algorithm to compute
the other homology groups was raised in \cite[Open Problem 1] {X20042}.
The present work is devoted to solving this problem.

In the direction of the more general problem of computing the homology 
of the category associated with a patial trace monoid action, there was carried out great work.
In \cite{X20083}, a complex for computing the Leech homology groups of trace monoids was obtained 
and the conjecture 
on the homological dimension of the augmented category \cite[Open Problem 2]{X20042} was confirmed.
Calculating the cohomological Leech dimension of trace monoids 
led to a generalization of Hilbert's Syzygy Theorem
to the polynomial rings in partially commuting variables \cite{X2011}.
For application, we tried to use the idea of \cite{win1995} to
add ``point at infinity'' in order to the action has become to be total.
This allowed us to build an algorithm for computing the homology groups 
of the augmented category of a partial trace monoid action on a set.
In \cite{X20082}, this result was used to find decomposition symptoms  of asynchronous systems.

But in general
these groups are not isomorphic to homology groups of the 
category associated with a partial trace monoid action. 
In addition, the unpleasant that even for a partial trace monoid action on single point, 
these groups can be
very large. In particular, they 
can contain an arbitrary finite Abelian subgroups in dimensions $n \geq 2$.
Haucourt \cite{hau2009} first to use the Baues-Wirsching cohomology groups 
 to study a mathematical model of concurrency.
Our solution indicates that even for
computing the integral homology groups of the category associated with a partial trace monoid action, it is  
helpful to use the Baues-Wirsching homology groups.

In the first section, we prove the auxiliary assertions. 
In the second, we study the Baues-Wirsching homology groups of 
a 
monoid action category and prove that they are reduced to the Leech homology groups. 
In the third, we get complexes to calculate the homology groups 
of the 
category assiciated with a 
trace monoid action category. 
The fourth part is devoted to algorithms and examples of calculating the
integer integral homology groups of a partial trace monoid action and 
 a CE net.

\tableofcontents

\section{Preliminaries}

In this section, we prove auxiliary Lemma \ref{mainlem} allowing us to reduce 
a
studying the Baues-Wirsching homology groups of the left fiber $h^{\mD}/X$ 
of Yoneda embedding
over $X \in \Set^{\mD^{op}}$ to the homology groups of the category $\mD$.

We deduce the formulas for finding the values of the functor obtained by the left
Kan extension from a functor defined 
on the category of factorizations of $h^{\mD}/X$
(Proposition \ref{desclan}).

\subsection{Notations}

Throughout the paper we use the following notation.
For any category $\mA$ and its objects $A, B \in \Ob\mA$, 
the 
set of all morphisms $A \to B$ is denoted by $\mA(A, B)$.
In the case of $\alpha \in \mA (A, B)$, we denote 
$\dom \alpha = A$ and $\cod \alpha = B $.
Let $\mA^{op}$ be the category opposite to $\mA$.
If $\mC$ is a small category, then $\mA^{\mC}$ denotes the category 
of functors $\mC\to \mA$ and natural transformations between them.
Let $\Set$ be the category of sets and maps (functions) and let $\Ab$ be the category 
of Abelian groups and homomorphisms.

Any monoid $M$ will be considered as a category with a single object.
A {\em right $M$-set $X$} or {\em right monoid action of a monoid $M$ on a set $X$} 
 is a functor which denoted by the same symbol $X: M^{op}\to \Set$.


\subsection{Left and right fibres of a functor}

Let ${\Phi}: \mC \to \mD $ be a functor from a small category $\mC$ to an arbitrary category $\mD$.
For any object $d \in \Ob \mD $ a {\em left fibre} ${\Phi}/ d$ \cite{gab1967}
(or {\em comma-category $\Phi \downarrow d$} \cite {mac1998}) is the category
whose objects are pairs $(c \in \Ob \mC, \alpha \in \mD ({\Phi} (c), d))$, and morphisms are
triples  
$(f \in \mC (c, c'), \alpha \in \mD ({\Phi} (c), d), \alpha' \in \mD ({\Phi} (c') , d)) $
making commutative the diagrams$$
\xymatrix{
{\Phi}(c) \ar[rd]_{{\Phi}(f)} \ar[rr]^{\alpha} && d\\
& {\Phi}(c') \ar[ru]_{\alpha'}
}
$$
Since in this case $\alpha=\alpha'{\Phi}(f)$, a morphism of a left fibre 
can be given as the pair
$(f,\alpha')$, but one should not forget the implicit existence of a morphism
 $\alpha=\alpha'{\Phi}(f)$.
The composition is defined as $(f',\alpha',\alpha'')\circ (f, \alpha, \alpha') =(f'f, \alpha, \alpha')$.
The identity morphism of an object 
 $(c,\alpha)$ equals $(1_c, \alpha, \alpha)$.

A {\em forgetful functor from a left fibre} $Q_d: {\Phi}/d\to \mC$ 
is defined as 
 assigning to
each $(c, \alpha)$ the object $c$ 
and to the morphism $(f, \alpha, \alpha')$ the morphism $f$.

A
{\em right fibre} $d/{\Phi}$ is opposite category to the left fibre ${\Phi}^{op}/d$ of a functor 
${\Phi}^{op}: \mC^{op}\to \mD^{op}$.  Its objects 
are given as pairs $(c, \alpha\in \mD( d,{\Phi}(c)))$,
and morphisms $(c,\alpha)\to (c',\alpha')$ are triples $(f,\alpha,\alpha')$
for which the following diagrams are commutative
$$
\xymatrix{
d \ar[rd]_{\alpha'} \ar[rr]^{\alpha} && {\Phi}(c) \ar[ld]^{{\Phi}(f)}\\
& {\Phi}(c')
}
$$

\subsection{Category of factorizations}

Let $\mC$ be a small category. A {\em category of factorizations $\fF\mC$}
is defined as follows. Its objects are morphisms of $\mC$, i.e. $\Ob\fF\mC= \Mor\mC$.
Every morphism $\alpha \to \beta$ in the category $\fF\mC$ is defined 
by a pair of morphisms
$(f: \dom\beta\to \dom\alpha, g:\cod\alpha\to \cod\beta)$, which makes commutative
diagram
$$
\xymatrix{
\cod\alpha \ar[r]^{g} & \cod\beta\\
\dom\alpha \ar[u]^{\alpha} & \dom\beta\ar[l]^f \ar[u]_{\beta}
}
$$ 
In this case, the morphism will be denoted by $\alpha\stackrel{(f,g)}\longrightarrow \beta$.
Composition of morphisms is defined by 
$(\beta\stackrel{(f',g')}\longrightarrow \gamma)\circ(\alpha\stackrel{(f,g)}\longrightarrow \beta)=
(\alpha\stackrel{(ff',g'g)}\longrightarrow \gamma$). Identity morphism equals
$1_{\alpha}=(\alpha\stackrel{(1_{\dom\alpha},1_{\cod\alpha})}\longrightarrow \alpha)$.

To each functor ${\Phi}: \mC\to \mD$ there corresponds a functor $\fF{{\Phi}}: \fF\mC \to \fF\mD$
given on objects as $\fF{{\Phi}}(\alpha)={\Phi}(\alpha)$ and defined 
on morphisms by the formula
$\fF{{\Phi}}(\alpha\stackrel{(f,g)}\to \beta)= 
({\Phi}(\alpha)\stackrel{({\Phi}(f),{\Phi}(g))}\longrightarrow {\Phi}(\beta))$.

\subsection{The forgetful functor from a left fibre of Yoneda embedding and its property}

Let $\mD$ be a small category and $X: \mD^{op}\to \Set$ a functor.
Consider the Yoneda embedding 
$$
	h^{\mD}: \mD\to \Set^{\mD^{op}}.
$$
Its values at $d\in \Ob\mD$ will be denoted by $h_d$.
It maps morphisms $\alpha: d_1\to d_2$ to the natural transformations 
$h_{\alpha}: h_{d_1}\to h_{d_2}$.
Since $X$ is the object of $\Set^{\mD^{op}}$,
there is the left fibre 
$h^{\mD}/X$ of  $h^{\mD}$ over $X$ and there is the forgetful functor 
$Q_X: h^{\mD}/X \to \mD$.
By Yoneda Lemma, there exists 
the bijection  $\Set^{\mD^{op}}(h_d, X)\stackrel{\cong}\to X(d)$, $\xi\mapsto \xi_d(1_d)$,
which is  natural in 
 $d\in \mD$ and  
$X\in \Set^{\mD^{op}}$. 
For $x\in X(d)$, denote by $\widetilde{x}: h_d\to X$ the natural transformation
corresponding to $x$.

A  category $\mC$ is {\em connected} if for any its objects $a,b\in \Ob\mC$, 
there exists 
a chain of morphisms in $\mC$ 
$$
	a\to b_1 \leftarrow a_1 \to b_2 \leftarrow \cdots \to b_n=b.
$$
{\em Connected component (of $a\in \Ob\mC$) in $\mC$} is a maximal connected subcategory
of $\mC$ (containing $a$). 

The following Lemma is crucial in the investigation of the Baues-Wirsching 
homology (and cohomology) groups for left fibres of the functor $h^{\mD}$.
Consider the functor between categories of factorizations 
$\fF{Q_X}: \fF({h^{\mD}/X})\to \fF{\mD}$ corresponding 
to the forgetful functor $Q_X: h^{\mD}/X\to \mD$.

\begin{lemma}\label{mainlem}
Let $\mD$ be a small category and let $X: \mD^{op}\to \Set$ be a functor.
Then for every $\alpha \in \Ob\fF\mD$,
each connected component
of the right fibre $\alpha/ \fF{Q_X}$ has an initial object.
\end{lemma}
{\sc Proof.} 
Objects in the category $\fF(h^{\mD}/X)$ are morphisms in $h^{\mD}/X$:
$$
\xymatrix{
	h_{d_1} \ar[rd]_{h_{\alpha}} \ar[rr]^{\widetilde{x}} && X\\
	&	h_{d_2}\ar[ru]_{\widetilde{y}} 
}
$$
Since $\widetilde{x}=\widetilde{y}h_{\alpha}$, 
we can describe
them with determining their pairs of morphisms 
$h_{d_1}\stackrel{h_{\alpha}}\to h_{d_2} \stackrel{\widetilde{y}}\to X$.
Morphisms in $\fF(h^{\mD}/X)$ are commutative diagrams
$$
\xymatrix{
	h_{d_2} \ar[rr]^{h_g} \ar[rd] && h_{e_2}\ar[ld]_{\widetilde{z}}\\
	 & X \\
	h_{d_1} \ar[uu]^{h_{\alpha}}\ar[ru] && h_{e_1}\ar[ll]^{h_f} \ar[uu]_{h_{\beta}} \ar[lu]
}
$$
In order to define these morphisms in $\fF(h^{\mD}/X)$,
it is enough to specify the morphisms
$h_{e_1}\stackrel{h_f}\to h_{d_1}\stackrel{h_{\alpha}}\to h_{d_2}\stackrel{h_g}\to
h_{e_2}\stackrel{\widetilde{z}}\to X$ in the category $\Set^{\mD^{op}}$.

In the category $\alpha/ \fF{Q_X}$, morphisms are pairs, each of which consists of an object 
$h_{d_1}\stackrel{h_{\alpha_1}}\to h_{d_2}\stackrel{\widetilde{z}}\to X$
in $\fF(h^{\mD}/X)$ and a morphism $\alpha \stackrel{(f,g)}\to 
\fF{Q_X}(h_{d_1}\stackrel{h_{\alpha_1}}\to h_{d_2}\stackrel{\widetilde{z}}\to X)$ категории $\fF\mD$.
They are defined by commutative diagrams

$$
\xymatrix{
e_2 \ar[r]^{g} & d_2\\
e_1 \ar[u]_{\alpha} & d_1 \ar[l]^f \ar[u]_{\alpha_1}
}
\qquad
\xymatrix{
h_{d_2} \ar[r]^{\widetilde{z}} & X\\
h_{d_1} \ar[u]^{h_{\alpha_1}} \ar[ru]_{\widetilde{z}h_{\alpha_1}}
}
$$
Taking this fact into account and using that the Yoneda embedding is full and $ \alpha_1 = g \alpha f $,
we can specify the objects of $ \alpha / \fF {Q_X}$ as the quadruples of morphisms
$h_{d_1}\stackrel{h_{f}}\to h_{e_1}\stackrel{h_{\alpha}}\to 
h_{e_2}\stackrel{h_{g}}\to h_{d_2}\stackrel{\widetilde{z}}\to X$.

In  $\alpha/\fF{Q_X}$, morphisms from the object 
$\alpha \stackrel{(f,g)}\to 
\fF{Q_X}(h_{d_1}\stackrel{h_{\alpha_1}}\to h_{d_2}\stackrel{\widetilde{z}}\to X)$
into the object 
$\alpha \stackrel{(f',g')}\to 
\fF{Q_X}(h_{d'_1}\stackrel{h_{\alpha'_1}}\to h_{d'_2}\stackrel{\widetilde{z'}}\to X)$
is specified as morphisms in  $\fF(h^{\mD}/X)$
$$
\xymatrix{
	h_{d_1}\ar[r]^{h_{\alpha_1}} & h_{d_2} \ar[r]^{\widetilde{z}} \ar[d]^{h_{\gamma}} & X\\
	h_{d'_1}\ar[r]^{h_{\alpha'_1}} \ar[u]^{h_{\beta}} & h_{d'_2} \ar[ru]_{\widetilde{z'}} 
}
$$
concerned with $(f,g)$ and $(f',g')$ by the commutative diagrams

$$
\xymatrix{
d_1 \ar[rd]^f \ar[rrr]^{\alpha_1} &&& d_2 \ar[dd]^{\gamma}\\
& e_1 \ar[r]^{\alpha} & e_2 \ar[ru]^g \ar[rd]_{g'}\\
d'_1 \ar[uu]^{\beta} \ar[ru]_{f'} \ar[rrr]_{\alpha'_1} &&& d_2
}
$$

It follows that the morphisms $\alpha/\fF{Q_X}$ 
can be specified by the commutative diagrams
$$
\xymatrix{
h_{d_1} \ar[r]^{h_f} & h_{e_1} \ar@{=}[d] \ar[r]^{h_{\alpha}} 
	& h_{e_2} \ar@{=}[d] \ar[r]^{h_g} &h_{d_2} \ar[d]^{h_{\gamma}} 
		\ar[r]^{\widetilde{z}} & X \ar@{=}[d]\\
h_{d'_1} \ar[u]^{h_{\beta}} \ar[r]^{h_{f'}} & h_{e_1} \ar[r]^{h_{\alpha}} 
	& h_{e_2} \ar[r]^{h_{g'}} &h_{d'_2} \ar[r]^{\widetilde{z'}} & X
}
$$
lines of which are quadruples of morphisms defining the objects of $\alpha/\fF{Q_X}$. 
If there is a morphism between objects of the category
$ \alpha / \fF {Q_X} $ specified by these quadruples, then 
$\widetilde{z}h_g=\widetilde{z'}h_{g'}$. 

Hence we conclude that the initial object of
connected component which contains the object
$h_{d_1}\stackrel{h_{f}}\to h_{e_1}\stackrel{h_{\alpha}}\to 
h_{e_2}\stackrel{h_{g}}\to h_{d_2}\stackrel{\widetilde{z}}\to X$
 in $ \alpha / \fF{Q_X}$,
 will be given by quadruple
$h_{e_1}\stackrel{=}\to h_{e_1}\stackrel{h_{\alpha}}\to 
h_{e_2}\stackrel{=}\to h_{e_2}\stackrel{\widetilde{z}h_g}\to X$.
Here, the equality symbol denotes the identity morphism.

Morphism from the initial object of the connected component into the object
$h_{d_1}\stackrel{h_{f}}\to h_{e_1}\stackrel{h_{\alpha}}\to 
h_{e_2}\stackrel{h_{g}}\to h_{d_2}\stackrel{\widetilde{z}}\to X$
will be given by the commutative diagram
$$
\xymatrix{
h_{e_1} \ar[r]^{=} & h_{e_1} \ar@{=}[d] \ar[r]^{h_{\alpha}} 
	& h_{e_2} \ar@{=}[d] \ar[r]^{=} &h_{e_2} \ar[d]^{h_{g}} 
		\ar[r]^{\widetilde{z}h_g} & X \ar@{=}[d]\\
h_{d_1} \ar[u]^{h_{f}} \ar[r]^{h_{f}} & h_{e_1} \ar[r]^{h_{\alpha}} 
	& h_{e_2} \ar[r]^{h_{g}} &h_{d_2} \ar[r]^{\widetilde{z}} & X
}
$$

\hfill $\Box$

\subsection{The left Kan extension along the forgetful functor of a left fibre}

For arbitrary functor $\Phi: \mC\to \mB$ and a cocomplete category $\mA$, 
a functor of the {\em left Kan extension} $Lan^{\Phi}: \mA^{\mC}\to \mA^{\mB}$ 
is defined as adjoint to the functor $\Phi^*: \mA^{\mB}\to \mA^{\mC}$ assigning 
to any functor $G: \mB\to \mA$ the composition $G\circ\Phi: \mC\to \mA$.
By \cite[Corollary X.3.4]{mac1998}, for every functor
$F: \mC\to \mA$, the values of  $Lan^{\Phi}F: \mB\to \mA$ can be described by formula
$$
(Lan^{\Phi}F) (b)= \coLim^{\Phi/b}(\Phi/b \stackrel{Q_b}\to \mC \stackrel{F}\to \mA)= \coLim^{\Phi/b}FQ_b.
$$
We will be working with the left Kan extension functor along the functor between 
the opposite categories
$\Phi^{op}: \mC^{op}\to \mB^{op}$. 
Since $\Phi^{op}/b \cong (b/\Phi)^{op}$, the left Kan extension of $F: \mC^{op}\to \Ab^{op}$
 takes the values
$(Lan^{\Phi^{op}}F) (b)= \coLim^{(b/\Phi)^{op}}FQ^{op}_b$.

Let $\mC$ be a small category whose every connected component 
has an initial object.  For $c \in \Ob \mC$, choose one initial object 
$i(c)$ in the connected component of  the category $\mC$,
 containing $c$. Denote the unique morphism from $i(c)$ to $c$ by
$i_c: i(c) \to c$.
Let $init(\mC)$ be the set of the chosen initial objects.

Let ${\Phi}: \mC\to \mB$ be a functor. For $b\in \Ob\mB$, objects 
$\zeta\in b/{\Phi}$ are specified by pairs $(c,\underline{\zeta})$, 
consisting of $c\in \Ob\mC$ with morphisms $\underline{\zeta}\in \mB(b,{\Phi}(c))$.

In some cases, it is convenient to omit the underscore symbol and write $\zeta$
instead of $\underline{\zeta}$. In particular for $\zeta=(c,\underline{\zeta})$  
and for arbitrary morphism $\sigma: a\to b$ категории $\mB$,
composition $\zeta\circ\sigma$ denotes the object $(c, \underline{\zeta}\circ\sigma)\in \Ob(a/{\Phi})$.

\begin{lemma}\label{compinit}\cite[Lemma 3.6]{X20081}
Let ${\Phi}: \mC\to \mB$ be a functor between small categories such 
that for every  $b\in \Ob\mB$ each connected component of  $b/{\Phi}$
has an initial object. Then for any functor $F: \mC^{op}\to \mA$,
its left Kan extension  $Lan^{{\Phi}^{op}}F: \mB \to \mA$ is isomorphic to the functor,
taking on   $b\in \Ob\mB$ the values $\coprod\limits_{\zeta\in init(b/{\Phi})}FQ_b^{op}(\zeta)$
and assigning to every morphism $\sigma: a\to b$ of $\mB$ the morphism  $\overline{\sigma}$
of $\mA$ determined for each $\zeta\in init(b/{\Phi})$ 
by the commutativity of the following diagrams: 
$$
\xymatrix{
\coprod\limits_{\zeta\in init(b/{\Phi})}FQ_b^{op}(\zeta) \ar[rr]^{\overline{\sigma}}
&& \coprod\limits_{\xi\in init(a/{\Phi})}FQ_a^{op}(\xi)\\
FQ_b^{op}(\zeta) = FQ_a^{op} (\zeta\circ\sigma)
\ar[rr]_{FQ_a^{op}(i_{\zeta\circ\sigma})} \ar[u]^{in_{\zeta}} && 
FQ_a^{op}(i(\zeta\circ\sigma)) \ar[u]_{in_{i(\zeta\circ\sigma)}}
}
$$
Here $in_{\zeta}$ are the canonical morphisms into the coproduct of objects.
\end{lemma}

\begin{proposition}\label{desclan}
Let $\mD$ be a small category and let $Q_X: h^{\mD}/X \to \mD$ be the left fibre forgetful functor.
For every functor $F: (\fF(h^{\mD}/X))^{op}\to \Ab$ 
and $\alpha\in \Ob\fF\mD$, we have 
$(Lan^{(\fF{Q_X})^{op}}F) (\alpha)= \bigoplus\limits_{x\in X(\cod\alpha)}F(\alpha,\widetilde{x})$.
The functor
$Lan^{(\fF{Q_X})^{op}}F$ 
assigns to each morphism $\alpha\stackrel{(f,g)}\longrightarrow \beta$
a morphism $(\overline{\alpha\stackrel{(f,g)}\longrightarrow \beta})$
that makes  the following diagram commutative for any
$y\in X(\cod\beta)$
$$
\xymatrix{
\bigoplus\limits_{y\in X(\cod\beta)}F(\beta,\widetilde{y}) 
\ar[rr]^{(\overline{\alpha\stackrel{(f,g)}\longrightarrow \beta})} && 
\bigoplus\limits_{x\in X(\cod\alpha)}F(\alpha,\widetilde{x})\\
F(\beta, \widetilde{y}) \ar[u]^{in_{(\beta,\widetilde{y})}} 
\ar[rr]_{ F(\alpha\stackrel{(f,g)}\longrightarrow \beta ,\widetilde{y}) } && 
F(\alpha, \widetilde{y}h_g) \ar[u]_{in_{(\alpha, \widetilde{y}h_g)}}
}
$$
\end{proposition}
{\sc Proof.} By Lemma \ref{mainlem}, connected components of 
all right fibres for the functor $\fF{Q_X}$, have initial objects and we can use 
Lemma \ref{compinit}.
For this aim, apply Lemma \ref{compinit} for $\mC=\fF(h^{\mD}/X)$, $\mB=\fF\mD$, $\mA=\Ab$.
Instead of ${\Phi}:\mC\to \mD$, we substitute the functor $\fF{Q_X}: \fF(h^{\mD}/X)\to \fF\mD$.
Substitute $b=\beta$, $a=\alpha$. 
We have obtained 
$Lan^{(\fF{Q_X})^{op}}F(\beta)= 
\bigoplus\limits_{\zeta\in init(\beta/\fF Q_X)}FQ^{op}_{\beta}(\zeta)$.
The objects of the category $\fF(h_*^{\mD}/X)$ defined by the diagrams 
$$
\xymatrix{
h_{d_2} \ar[r]^{\widetilde{y}} & X\\
h_{d_1}\ar[u]^{h_{\beta}} \ar[ru]_{\widetilde{y}h_{\beta}}
}
$$
will be given as the pairs $(\beta,\widetilde{y})$.

Consider the right fibre forgetful functor
$Q_{\beta}: \beta/\fF{Q_X}\to \fF(h_*^{\mD}/X)$.
For $\zeta \in \beta/\fF{Q_X}$ given by four natural transformations
$$
h_{d'_1}\stackrel{h_{f'}}\longrightarrow
h_{e_1}\stackrel{h_{\beta}}\longrightarrow h_{e_2}\stackrel{h_{g'}}\longrightarrow
h_{d'_2}\stackrel{\widetilde{y'}}\longrightarrow X
$$
we have $Q_{\beta}(\zeta)=(g'\beta f', \widetilde{y})$.
If $\zeta\in init(\beta/\fF Q_X)$, then 
$\zeta$ is given by the quadruple 
\begin{equation}\label{fourz}
h_{\dom{\beta}}\stackrel{=}\longrightarrow
h_{\dom{\beta}}\stackrel{h_{\beta}}\longrightarrow h_{\cod\beta}\stackrel{=}\longrightarrow
h_{\cod{\beta}}\stackrel{\widetilde{y}}\longrightarrow X
\end{equation}
Hence, for fixed $\beta$, each object $\zeta\in init(\beta/\fF Q_X)$ is determined by the morphism
$\widetilde{y}: h_{\cod\beta}\to X$,  or equivalently, by $y\in X(\cod\beta)$.
Therefore,
$$
Lan^{(\fF{Q_X})^{op}}F(\beta)= 
\bigoplus\limits_{y\in X(\cod\beta)}F(\beta, \widetilde{y}).
$$
Consider a morphism
$\sigma=(\alpha\stackrel{(f,g)}{\longrightarrow}\beta)$ in the category $\fF\mD$
and describe the map $Lan^{(\fF{Q_X})^{op}}F (\sigma)$.
For $\zeta$, given by the quadruple (\ref{fourz}), the composition $\zeta\circ \sigma$ is equal to
$$
h_{d_1}\stackrel{h_{f}}\longrightarrow
h_{e_1}\stackrel{h_{\beta}}\longrightarrow h_{e_2}\stackrel{h_{g}}\longrightarrow
h_{d_2}\stackrel{\widetilde{y}}\longrightarrow X.
$$
It follows from $g\alpha f=\beta$ that 
$FQ^{op}_{\alpha}(\zeta\circ\sigma)=FQ^{op}_{\beta}(\zeta)=F(\beta,\widetilde{y})$.
The morphism $i_{\zeta\circ\sigma}: i(\zeta\circ\sigma)\to \zeta\circ\sigma$ is given
by the commutative diagram
$$
\xymatrix{
h_{\dom\alpha} \ar[r]^{=} & h_{\dom\alpha} \ar@{=}[d] \ar[r]^{h_{\alpha}} 
	& h_{\cod\alpha} \ar@{=}[d] \ar[r]^{=} &h_{\cod\alpha} \ar[d]^{h_{g}} 
		\ar[r]^{\widetilde{y}h_g} & X \ar@{=}[d]\\
h_{\dom\beta} \ar[u]^{h_{f}} \ar[r]^{h_{f}} & h_{\dom\alpha} \ar[r]^{h_{\alpha}} 
	& h_{\cod\alpha} \ar[r]^{h_{g}} &h_{\cod\beta} \ar[r]^{\widetilde{y}} & X
}
$$
Thus, $FQ^{op}_{\alpha}(i_{\zeta\circ\sigma}): FQ^{op}_{\alpha}(\zeta\circ\sigma)\to 
FQ^{op}_{\alpha}(i(\zeta\circ\sigma))$ is equal to 
$F(\alpha\stackrel{(f,g)}\longrightarrow\beta): F(\beta,\widetilde{y})\to  F(\alpha,\widetilde{y}h_g)$.
\hfill $\Box$

\section{Baues-Wirshing homology of a monoid action category}

First, we define the Baues-Wirsching homology groups and observe the required properties.
Considering a monoid $M$ as a category with a single object, and a right $M$-set as
a functor $X: M^{op} \to \Set$, we obtain the category $h^{M}_*/X$.
A {\em category $K(X)$ associated with a monoid action of $M$ on $X$} or {\em a monoid action category $K(X)$}
can be defined as the opposite category to $h^{M}_*/X$.
We will prove that for any convex subcategory of the monoid action category,
a computing the Baues-Wirsching homology can be reduced to calculating 
 the Leech homology groups of the monoid $M$.
And then we apply this statement to the homology of category associated 
with a partial monoid action.

\subsection{Homology of small categories with coefficients in functors}

Given a small category $\mC$, we consider the category $\Ab^{\mC^{op}}$ of functors  
$F: \mC^{op}\to \Ab$. It is well known that this category is Abelian that has 
enough projectives (and enough injectives). 
Hence, for any right exact functor $\Ab^{\mC^{op}}\to \Ab$, 
there are its left derived functors.
This category has exact coproducts (and exact products).

Define the homology groups of the category $\mC$ with coefficients in $F$ 
by left derived functors 
$\coLim_n^{\mC^{op}}: \Ab^{\mC^{op}}\to \Ab$ of the colimit functor
 $\coLim^{\mC^{op}}: \Ab^{\mC^{op}}\to \Ab$.

\begin{definition}
For $n\geq 0$, Abelian groups 
$\coLim_n^{\mC^{op}}F$ are called $n$th homology groups of $\mC$ 
with coefficients in $F$.
\end{definition}

Formally, for a functor $F$ considered as the object in $\Ab^{\mC^{op}}$, 
we take an arbitrary projective resolution
$$
0 \leftarrow F \leftarrow F_0 \leftarrow F_1 \leftarrow \cdots \leftarrow F_n \leftarrow \cdots
$$
in the category $\Ab^{\mC^{op}}$ 
and groups $H_n(\mC^{op},F)$ are defined as the factor groups $Ker d_n/ Im{d_{n+1}}$ of the complex
obtained by the application of the colimit functor
$$
0 \stackrel{d_0}\leftarrow \coLim^{\mC^{op}}F_0 \stackrel{d_1}\leftarrow \coLim^{\mC^{op}}F_1 
\stackrel{d_2}\leftarrow \cdots 
\stackrel{d_n}\leftarrow \coLim^{\mC^{op}}F_n \leftarrow \cdots.
$$

\begin{proposition}\label{liftcol}
Let ${\Phi}: \mC\to \mB$ be a functor between small categories such that
for every $b\in \Ob\mB$, each connected component of $b/{\Phi}$
has an initial object. Then, for any functor $F: \mC^{op}\to \Ab$,
there are isomorphisms $\coLim_n^{\mC^{op}}F\cong \coLim_n^{\mB^{op}}Lan^{{\Phi}^{op}}F$.
\end{proposition}
{\sc Proof.}  Take any projective resolution of the functor $F\in \Ab^{\mC^{op}}$
in $\Ab^{\mC^{op}}$. Since $Lan^{{\Phi}^{op}}$ is a left adjoint to the exact functor
$\Ab^{\mB^{op}}\to \Ab^{\mC^{op}}$ acting as $G\mapsto G\circ {\Phi}^{op}$, it takes 
projective objects to projective objects. The functor $Lan^{{\Phi}^{op}}$ translates a projective 
resolution $0\leftarrow F \leftarrow F_0 \leftarrow F_1 \leftarrow \cdots$ to the complex
\begin{equation}\label{exlan}
0 \leftarrow Lan^{{\Phi}^{op}}F \leftarrow Lan^{{\Phi}^{op}}F_0 \leftarrow \cdots  
\leftarrow Lan^{{\Phi}^{op}}F_n \leftarrow \cdots
\end{equation}
the members of which have values
$Lan^{{\Phi}^{op}}F(b)= \bigoplus\limits_{\zeta\in init(b/{\Phi})} FQ^{op}_b (\zeta)$.
Hence, for every $b\in \Ob\mB$, the sequence of values of the members at $b$ is exact.
Therefore, the sequence of functors (\ref{exlan}) is exact.
By applying the functor $\coLim^{\mB^{op}}$ to the resulting resolution,
we obtain a complex 
whose homology groups are isomorphic to $\coLim_n^{\mB^{op}}Lan^{{\Phi}^{op}} F$. 
Since $\coLim^{\mB^{op}}Lan^{{\Phi}^{op}} F\cong \coLim^{\mC^{op}}F$, the 
homology groups of the obtained complex are isomorphic to $\coLim_n^{\mC^{op}} F$. 
Consequently,  
$\coLim_n^{\mC^{op}} F\cong \coLim_n^{\mB^{op}}Lan^{{\Phi}^{op}} F$.
\hfill $\Box$

Let $\mC$ be a small category. A full subcategory $\mD\subseteq \mC$ is  
{\em closed} if for every morphism $c\to d$ from $c\in \Ob\mC$ into $d\in \Ob\mD$
will be true that $c\in \Ob\mD$. In this case, since $\mD$ is full, the morphism
 $c\to d$ will be belong to $\mD$.

Let ${\Phi}: \mD\stackrel{\subseteq}\to \mC$ be the inclusion of a closed subcategory. 
Consider a functor $F: \mD^{op}\to \Ab$.
Extending the values of $F$ 
on objects and morphisms by zeros, we obtain the functor $F_{\mD}: \mC^{op} \to \Ab $,
which is called the {\em functor obtained by adding zeros}.


For each $c \in \Ob \mC$, the right fibre $c /{\Phi} $ is either empty
(and hence $ Lan^{{\Phi} ^ {op}} F (c) = \coLim^ {(c / {\Phi}) ^ {op}} FQ^{op}_c = 0 $)
or has an initial object $ (c, 1_c) $ 
(and $ Lan^{{\Phi} ^ {op}} F (c) = \coLim ^ {(c / {\Phi})^{op}} FQ^{op} _c = F (c) $).
Thus $F_{\mC}=Lan^{{\Phi}^{op}}{F}$. 
It follows from Proposition \ref{liftcol} the following
\begin{corollary}
If $\mD\subseteq \mC$ is a closed subcategory, then
$\coLim_n^{\mC^{op}}F_{\mD}\cong \coLim_n^{\mD^{op}}F$ for any functor $F: \mD^{op}\to \Ab$.
\end{corollary}


\begin{definition}
Let $\mC$ be a small category and let $F: \fF\mC^{op}\to \Ab$ be a functor. 
{\em Baues-Wirsching homology groups of $\mC$ with coefficients in $F$} are the Abelian 
groups $\coLim_n^{\fF\mC^{op}}F$. In particular, 
for an arbitrary monoid $M$ considered as a category with unique object,
Baues-Wirsching homology groups with coefficients in $F: \fF{M}^{op}\to Ab$
are called the Leech homology groups of $M$.
\end{definition}

Full subcategory $\mD\subseteq \mC$ is called  {\em convex}
if, together with any of its objects $d, d'\in \Ob\mD $, it contains
every object $c \in \Ob\mC$ for which there exists a pair of morphisms $ d \to c \to d'$ in the category
$\mC$.
In this case the subcategory
$\fF\mD\subseteq \fF\mC$ is closed.
This implies 

\begin{corollary}\label{liftconvex}
If $\mD\subseteq \mC$ is a convex subcategory, then
$\coLim_n^{\fF\mC^{op}}F_{\fF\mD}\cong \coLim_n^{\fF\mD^{op}}F$ for any functor
$F: \fF\mD^{op}\to \Ab$.
\end{corollary} 

\subsection{The homology groups of a (total) monoid action category}

For any monoid $M$, the category $h_*^{M}/X$ will be denoted by $M/X$.

For $X\in \Set^{M^{op}}$, 
a {\em right monoid action category $K(X)$} of $M$ on $X$ is a small category with the set 
of objects $\Ob{K(X)}=X$. Its morphisms $x\to y$ are defined as triples $(\mu,x,y)$, 
$\mu\in M$, satisfying to $x\cdot\mu= y$. Composition is defined by the formulas
$(\nu, y, z)\circ (\mu, x, y)= (\mu\nu, x, z)$ for $x,y,z\in X$, $\mu, \nu\in M$.
We will consider a morphism $(\mu, x, y)$ as the pair $(x,\mu)$. The composition 
will be defined as
$(x\cdot\mu,\nu)\circ(x,\mu)= (x, \mu\nu)$. The objects in the category $M/X$ are 
natural transformation $\widetilde{x}: h_{M}\to X$, and morphisms are 
commutative triangles
$$
\xymatrix{
	h_{M} \ar[rd]_{h_{\mu}} \ar[rr]^{\widetilde{x}} && X\\
	& 	h_{M} \ar[ru]_{\widetilde{y}}
}
$$
determined by pairs $(\mu,\widetilde{y})$. This follows that the functor acting
$\widetilde{y}\mapsto y$ on objects and $(\mu,\widetilde{y})\mapsto (y,\mu)$ on morphisms,
is an isomorphism of categories $(M/X)^{op} \stackrel{\cong}\to K(X)$.

The category $K(X)$ is opposite to $M/X$, but the category of factorizations are isomorphic.
We construct the isomorphism.

Morphisms $(\mu, x, x\cdot\mu)$ of $K(X)$ are given by pairs $(x,\mu)$.
Hence, the objects of $\fF K(X)$ are given by pairs $(x,\mu)$ where  $x\in X$, $\mu\in M$.  
Morphisms $(x,\mu)\stackrel{(g,f)}\to (y,\nu)$ are commutative diagrams in $K(X)$
$$
\xymatrix{
x\cdot\mu \ar[r]^{f} & y\cdot \nu\\
x \ar[u]^{\mu} & y \ar[l]^{g}\ar[u]_{\nu}
}
$$


\begin{proposition}\label{ptheta}
Isomorphism $\Theta: \fF(M/X)\stackrel{\cong}\to \fF{K(X)}$
is realized by
the functor that maps objects by the formula 
$ \Theta (\mu, \widetilde {x}) = (x, \mu)$, 
and acting on morphisms as
$$
\Theta(\mu\stackrel{(f,g)}\longrightarrow \nu, \widetilde{y})=
((y\cdot g, \mu) \stackrel{(g,f)}\longrightarrow (y,\nu)).
$$
\end{proposition}
{\sc Proof.} 
The commutative diagram in $K(X)$
$$
\xymatrix{
y\cdot g\mu \ar[r]^f & y\cdot \nu\\
y\cdot g \ar[u]^{\mu} & y \ar[l]^g \ar[u]_{\nu}
}
$$
clearly shows that $(y\cdot g, \mu) \stackrel{(g,f)}\longrightarrow (y,\nu)$  
defines the morphism in the category $\fF{K(X)}$.
\hfill $\Box$

Let $X$ be a right $M$-set and let $S\subseteq X$ be an arbitrary subset.
 Denote by $K(S)\subseteq K(X)$
the full subcategory with the set of objects $S$.

The following proposition shows
that the Baues-Wirsching homology groups of a monoid action category  can be
expressed in terms of the Leech homology groups.

Above, we considered the functors
$$
\xymatrix{
&	\fF(M/X)^{op}\ar[d]^{\Theta^{op}} \ar[rd]^{\fF{Q_X}^{op}}\\
\fF{K(S)}^{op} \ar[r]^{\subseteq} & \fF{K(X)^{op}} & \fF{M}^{op}
}
$$
\begin{proposition}\label{bauesleech}
Let $M$ be an arbitrary monoid with right action on a set $X$.
If $K(S)\subseteq K(X)$ is a convex subcategory and 
$G: \fF{K(S)}^{op} \to \Ab$ a functor, then there are isomorphisms
$$
\coLim_n^{\fF{K(S)}^{op}}G \cong \coLim_n^{\fF{M}^{op}}Lan^{\fF{Q_X}^{op}}
(G_{\fF{K(S)}}\circ\Theta^{op}),
$$
for all $n\geq 0$.
\end{proposition} 
{\sc Proof.} 
This follows from the isomorphisms
\begin{multline*}
\coLim_n^{\fF{K(S)}^{op}}G \cong \coLim_n^{\fF{K(X)}^{op}}G_{K(S)} \mbox{ (Corollary \ref{liftconvex})}\\
	\cong \coLim_n^{\fF{(M/X)}^{op}}(G_{K(S)}\circ\Theta^{op}) \mbox{ (Proposition \ref{ptheta})}\\
	\cong \coLim_n^{\fF{M}^{op}}Lan^{\fF{Q_X}^{op}}
(G_{\fF{K(S)}}\circ\Theta^{op}) \\
\mbox{ (Proposition \ref{liftcol} and Lemma \ref{mainlem})}
\end{multline*}
\hfill $\Box$

\subsection{The homology groups of a partial monoid action category}

A {\em partial map} of sets $ A\stackrel{f}\rightharpoonup B$ is given 
by a relation  $f\subseteq A\times B$ such that 
$$
	(x,y_1)\in f ~\&~ (x,y_2)\in f ~\Rightarrow~ y_1=y_2.
$$
The composition $(A\stackrel{g\circ f}\rightharpoonup C)$= 
$(B \stackrel{g}\rightharpoonup C)\circ(A\stackrel{f}\rightharpoonup B)$ 
is defined by the composition of relations
$$
	g\circ f= \{(x,z)\in A\times C| (\exists y\in B) (x,y)\in f \& (y,z)\in g\}.
$$
The identity morphism $1_A: A\to A$ equals $\{(x,x)| x\in A\}$.
Let $P\Set$ be the category of sets and partial maps.

A  {\em partial monoid action of monoid $M$} on a set $S$ is a functor $M^{op}\to P\Set$.
Following \cite{win1995}, consider the category  $\Set_*$ of sets 
$S_*=S\sqcup \{*\}$ with  the point $*$ and  
maps $f: S_*\to S_*$ satisfying 
$f(*)=*$. The point $*$ is identical for all objects.
There is an isomorphism
$P\Set\stackrel{\cong}\to \Set_*$. This follows that  a partial monoid action of $M$ on $S$
can be considered as the functor $S_*: M^{op}\to \Set_*$.
Denote $x\cdot\mu=S_*(\mu)(x)$. The operation '$\cdot$' has the following properties:
\begin{enumerate}
\item $(\forall x\in S_*)~ x\cdot 1=x$,
\item $(\forall x\in S_*)(\forall \mu, \nu \in M)~ x\cdot(\mu\nu)= (x\cdot\mu)\cdot\nu$,
\item \label{thirdproper} $(\forall \mu\in M)~ *\cdot \mu=*$.
\end{enumerate}

Let $U: \Set_*\to \Set$ be a functor defined by $U(S_*)=S_*$ and $U(f)=f$.
For a monoid $M$ and a functor $S_*: M^{op}\to \Set_*$, denote by $US_*: M^{op}\to \Set$
 the composition $U\circ S_*$. Denote by $K_*(S)$ the category $K(US_*)$.

The property (\ref{thirdproper}) leads to the implication:
$x\cdot\mu\not=* \Rightarrow x\not=*$.
This implies 

\begin{lemma}\label{KSclosed}
The subcategory $K(S)\subseteq K_*(S)$ is closed.
\end{lemma}

A closed subcategory is convex. 
By substituting in Proposition \ref{bauesleech}, for $X=US_*$, we obtain
the following

\begin{corollary}
For any  partial monoid action $S_*: M^{op}\to \Set_*$
and a functor
$G: \fF{K(S)}^{op} \to \Ab$  
there are isomorphisms
$$
\coLim_n^{\fF{K(S)}^{op}}G \cong \coLim_n^{\fF{M}^{op}}Lan^{\fF{Q_{US_*}}^{op}}
(G_{\fF{K(S)}}\circ\Theta^{op}),
$$
for all $n\geq 0$.
\end{corollary}

\section{The homology groups of a trace monoid action category}

We turn to our main results. In the first part, we will construct
complex to compute the Baues-Wirsching homology groups of the category
$\fF(M(E, I)/X)$ using the Leech homology groups
of the monoid $M(E,I)$.
In the second, study the homology of the convex subcategories
of the trace monoid action category.

\subsection{Baues-Wirsching homology of a left fibre of Yoneda embedding}

Recall the definition of a trace monoid.
Let $E$ be a set which, unlike in \cite{die1997}, we do not assume to be finite.
Subset $I \subseteq E \times E$ is
called {\em independence relation} if it is irreflexive and symmetric, i.e.
the following conditions are satisfied:

\begin{itemize}
\item $(\forall a\in E)~ (a,a)\notin I$,
\item $(\forall a\in E)(\forall b\in E)~ (a,b)\in I \Rightarrow (b,a)\in I$.
\end{itemize}

Elements $a,b\in E$ are  {\em independent} if  $(a,b)\in I$. 
Let $E^*$ be the free monoid of words
 including the empty word $1$ with the concatenation operation 
$(v,w)\mapsto vw$.
For an arbitrary independence relation $I\subseteq E\times E$,
define the equivalence relation $\equiv_I$ on $E^*$ 
for which $w_1\equiv_I w_2$ if the word $w_2$ can be obtained from $w_1$ by 
a finite sequence of 
permutations of adjacent independent elements.
For any $w\in E^*$, its equivalence class $[w]=\{v\in E^*| v\equiv_I w\}$  
 consists of a finite number of words. In particular, the class of the empty word 
$1$ equals $\{1\}$.
A {\em trace monoid $M(E,I)$} is the factor set $E^*/(\equiv_I)$
with the operation $[v][w]=[vw]$. 
Its neutral element $ \{1\} $ denoted by $1$.

A trace monoid $M(E,I)$ is {\em locally finite dimensional} if $E$ does not contain 
infinite subsets of pairwise independent elements.
Consider an arbitrary total order relation on $E$.
For $n\geq 1$, denote 
$$
T_n(E,I)=\{(a_1, \cdots, a_n)| a_i<a_j \mbox{ and }
(a_i,a_j)\in I, \mbox{ for all } 1\leq i< j\leq n\}.
$$
Let $T_0(E,I)= \{1\}$.

\begin{theorem}\label{main1}
Let $M(E,I)$ be a locally finite dimensional trace monoid 
and let 
$X:M(E,I)^{op}\to \Set$ be a right action of $M(E,I)$ on a set $X$.
Then, for any functor $F: (\fF(M(E,I)/X))^{op}\to \Ab$, the groups
 $\coLim_n^{(\fF(M(E,I)/X))^{op}} F$ will be isomorphic to the homology groups 
of the complex
\begin{multline*}
0 \leftarrow \bigoplus\limits_{x\in X} F(1, \widetilde{x}) \stackrel{d_1}\longleftarrow
\bigoplus\limits_{(x,a_1)\in X\times T_1(E,I)} F(a_1, \widetilde{x}) \leftarrow \cdots\\
\cdots \leftarrow 
\bigoplus\limits_{(x,a_1, \cdots, a_{n-1})\in X\times T_{n-1}(E,I)} 
F(a_1\cdots a_{n-1}, \widetilde{x})\\
 \stackrel{d_n}\longleftarrow
\bigoplus\limits_{(x,a_1, \cdots, a_{n})\in X\times T_{n}(E,I)} 
F(a_1\cdots a_{n}, \widetilde{x})\leftarrow \cdots
\end{multline*}
with the differentials $d_n$ 
defined on direct summands as
\begin{multline*}
	d_n(x, a_1, \cdots, a_n, \psi)= \\
		\sum_{i=1}^n (-1)^i (x a_i, a_1, \cdots, a_{i-1}, a_{i+1}, \cdots, a_n,
			F((\widehat{1,a_i}), 
				\widetilde{x})(\psi) )\\
	-\sum_{i=1}^n (-1)^i (x, a_1, \cdots, a_{i-1}, a_{i+1}, \cdots, a_n,
			F((\widehat{a_i,1}), 
				\widetilde{x})(\psi) )
\end{multline*}
for all $\psi\in F(a_1\cdots a_n, \widetilde{x})$.
Here are used the notations 
\begin{gather*}
(\widehat{1,a_i})= 
(a_1 \cdots a_{i-1} a_{i+1} \cdots a_n\stackrel{(1,a_i)}\to a_1\cdots a_n),\\ 
(\widehat{a_i,1})=
(a_1 \cdots a_{i-1} a_{i+1} \cdots a_n\stackrel{(a_i,1)}\to a_1\cdots a_n).
\end{gather*}
\end{theorem}
{\sc Proof.}
In order to find the homology groups of the category $\fF (M (E, I) / X)$
with coefficients in $ F: \fF (M (E, I) / X)^{op} \to \Ab$,
consider the left Kan extension $Lan^{(\fF {Q_X)}^{op}} F$.
By Lemma \ref{mainlem}, for $\alpha\in \Ob\fF(M(E,I))$, 
every connected component of the right fibre $\alpha/\fF{Q_X}$ 
has an initial object. Hence, the functor $\fF{Q_X}$ satisfies to conditions of Proposition
 \ref{liftcol}. Applying Proposition \ref{liftcol}, we obtain an isomorphism

$$
\coLim^{\fF(M(E,I)/X)^{op}}_n F\cong \coLim^{\fF{M(E,I)}^{op}}_n Lan^{\fF{Q_X}^{op}} F.
$$

In accordance with \cite[Theorem 2.16]{X20084}, for any functor
$G: \fF {M (E, I)}^{op} \to \Ab $, Abelian groups
$ \coLim^{\ fF {M (E, I)} ^ {op}} _n G $
are isomorphic to 
homology groups of the complex

\begin{multline*}
0 \leftarrow  G(1) \stackrel{d_1}\longleftarrow
\bigoplus\limits_{a_1\in T_1(E,I)} G(a_1) \leftarrow \cdots\\
\cdots \leftarrow 
\bigoplus\limits_{(a_1, \cdots, a_{n-1})\in T_{n-1}(E,I)} 
G(a_1\cdots a_{n-1})
 \stackrel{d_n}\longleftarrow
\bigoplus\limits_{(a_1, \cdots, a_{n})\in T_{n}(E,I)} 
G(a_1\cdots a_{n})\leftarrow \cdots
\end{multline*}
with differentials
\begin{multline}\label{diffref}
d_n(a_1, \cdots, a_n, \varphi)= \\
\sum\limits_{i=1}^n (-1)^i(a_1,\cdots, a_{i-1},a_{i+1},\cdots, a_n,
	G(\widehat{1,a_i})(\varphi)-G(\widehat{a_i,1})(\varphi))
\end{multline}

Substitute $G= Lan^{\fF{Q_X}^{op}} F$. Since $M(E,I)$ is the category with a single object,
the functor $X$ has a unique value on the objects. Denote this value by $X$.
From Proposition  \ref{desclan} we obtain
 $G(1)=\bigoplus\limits_{x\in X}F(1,\widetilde{x})$, $\cdots$, 
$G(a_1 \cdots a_n)= \bigoplus\limits_{x\in X}F(a_1 \cdots a_n,\widetilde{x})$, $\cdots$.
We have arrived to the complex consisting of the Abelian groups
$$
C_n= \bigoplus\limits_{(x,a_1, \cdots, a_{n})\in X\times T_{n}(E,I)} 
F(a_1\cdots a_{n}, \widetilde{x})
$$
Let us describe the differentials $d_n: C_n\to C_{n-1}$.

To this aim we note that 
for any small category $\mD$, a functor $X:\mD^{op} \to \Set$, objects $d,d'\in \Ob \mD$,
and a morphism $ g \in \mD (d, d')$,  the equality $\widetilde{y}h_{g}=\widetilde{X(g)(y)}$ holds.
Hence, when $ \mD = M(E, I)$ is a trace monoid with a right action on  $X$,
we have $ \widetilde {y} h_{g} = \widetilde {yg}$ for all $g \in M(E, I)$.

We make also remark on the designation of the elements of direct summands.
For $i\in J$, the elements of a direct summand $A_i$ in $\bigoplus\limits_{j\in J}A_j$ 
are denoted by $(i,a)$.
So, the inclusion  $in_i: A_i\to \bigoplus\limits_{j\in J}A_j$ is defined as $in_i(a)=(i,a)$.
In particular, the direct summand in
$\bigoplus\limits_{x\in X}F(a_1\cdots a_n, \widetilde{x})$ corresponding to index
$y\in X$ consists of pairs $(y,\psi)$ where $\psi\in F(a_1\cdots a_n)$.
A direct summand in 
$\bigoplus\limits_{{(x,a_1, \cdots, a_n)}\in X\times T_n(E,I)} F(a_1\cdots a_n,\widetilde{x})$
corresponding to index $(a_1, \cdots, a_n)$ consists of elements 
$(a_1, \cdots, a_n, \varphi)$ where $\varphi\in \bigoplus\limits_{x\in X}F(a_1\cdots a_n,\widetilde{x})$.

Formulas (\ref{diffref}) describe differentials $d_n: C_n\to C_{n-1}$ 
 by homomorphisms of the functor $G$. We must express $d_n$ by homomorphisms 
of $F$.
Since $G$ is the left Kan extension of  $F$, we can use Proposition \ref{desclan}.
Consider an arbitrary $\varphi$ from a direct summand of
$\bigoplus\limits_{x\in X}F(a_1\cdots a_n, \widetilde{x})$ and describe how to operate the 
homomorphisms 
$G(\widehat{1,a_i})$ and $G(\widehat{a_i,1})$ on it. By Proposition \ref{desclan}, 
for $\varphi=(y,\psi)$ with $\psi\in F(a_1\cdots a_n, \widetilde{y})$,
the following formulas hold
\begin{gather}
\label{form1} G(\widehat{1,a_i})(y,\psi)= (ya_i, F((\widehat{1,a_i}), \widetilde{y})(\psi))\\
 G(\widehat{a_i,1})(y,\psi)= (y, F((\widehat{a_i,1}), \widetilde{y})(\psi))
\label{form2}
\end{gather}
It is clear that
\begin{multline*}
d_n(a_1, \cdots, a_n, \varphi)= \\
\sum\limits_{i=1}^n (-1)^i(a_1,\cdots, a_{i-1},a_{i+1},\cdots, a_n,
	G(\widehat{1,a_i})(\varphi))\\
- \sum\limits_{i=1}^n (-1)^i(a_1,\cdots, a_{i-1},a_{i+1},\cdots, a_n, G(\widehat{a_i,1})(\varphi))
\end{multline*}
Substituting  $\varphi=(y,\psi)$ and using formulas (\ref{form1})-(\ref{form2}) complete the proof. 
\hfill $\Box$

\subsection{Baues-Wirsching homology of convex subcategories in the trace monoid action category}

Given a right trace monoid action of  $M(E,I)$ on $X$, 
consider a subset $S\subseteq X$.
Morphisms in the category $K(S)$ are exhausted pairs $(x,\mu)$ such that 
$\mu\in M(E,I)$, $x\in S$, and $x\cdot\mu\in S$.

If $K(S)\subseteq K(X)$ is a convex subcategory, then for an arbitrary total order on $E$, 
we can define the sets
$$
Q_n(S,E,I)=\{(x,a_1, \cdots, a_n)\in S\times T_n(E,I)| x\cdot a_1\cdots a_n\in S\} 
\mbox{ for all } n\geq 0.
$$
In particular, $Q_0(S,E,I)=S$.

\begin{theorem}\label{solution}
Let $M(E,I)$ be a locally finite dimensional trace monoid
and let $X: M(E,I)^{op}\to \Set$ be a right action. Let
 $S\subseteq X$ be a subset for which  $K(S)\subseteq K(X)$ is a convex subcategory.
Then, for an arbitrary total order on  $E$ and any functor $G: \fF{K(S)}^{op}\to \Ab$, 
Abelian groups $\coLim_n^{\fF{K(S)}^{op}}G$ will be isomorphic to homology groups 
of the complex
\begin{multline}\label{maincomp}
0 \leftarrow \bigoplus\limits_{x\in S} G(x, 1) \stackrel{d_1}\longleftarrow
\bigoplus\limits_{(x,a_1)\in Q_1(S,E,I)} G(x, a_1) \leftarrow \cdots\\
\cdots \leftarrow 
\bigoplus\limits_{(x,a_1, \cdots, a_{n-1})\in Q_{n-1}(S, E,I)} 
G(x, a_1\cdots a_{n-1})\\
 \stackrel{d_n}\longleftarrow
\bigoplus\limits_{(x,a_1, \cdots, a_{n})\in  Q_{n}(S,E,I)} 
G(x, a_1\cdots a_{n})\leftarrow \cdots 
\end{multline}
with defferentials $d_n$ 
defined on the direct summands by
\begin{multline}\label{maindiff}
	d_n(x, a_1, \cdots, a_n, \psi)= \\
		\sum_{i=1}^n (-1)^i (x\cdot a_i, a_1, \cdots, a_{i-1}, a_{i+1}, \cdots, a_n,\\
			G((x\cdot a_i, a_1 \cdots a_{i-1} a_{i+1} \cdots a_n) \stackrel{(a_i,1)}\to
(x, a_1 \cdots a_n))(\psi))\\
	-\sum_{i=1}^n (-1)^i (x, a_1, \cdots, a_{i-1}, a_{i+1}, \cdots, a_n,\\
			G((x, a_1 \cdots a_{i-1} a_{i+1} \cdots a_n) \stackrel{(1,a_i)}\to
(x, a_1 \cdots a_n))(\psi))\\
\end{multline}
for all $\psi\in G(x, a_1\cdots a_n)$.
\end{theorem}
{\sc Proof.} By Proposition \ref{bauesleech},
we have isomorphisms
$$
\coLim_n^{\fF{K(S)}^{op}}G \cong \coLim_n^{\fF{M(E,I)}^{op}}Lan^{\fF{Q_X}^{op}}
(G_{\fF{K(S)}}\Theta^{op})
$$ 
for all $n\geq 0$.
Application of Theorem \ref{main1} to the functor
$F=G_{\fF{K(S)}}\circ\Theta^{op}$
leads us to the complex (\ref{maincomp}).

We will prove formulas (\ref{maindiff}).
The morphism $\Theta((\widehat{1,a_i}), \widetilde{x})$ equals to the morphism  
$(x\cdot a_i, a_1\cdots a_{i-1}a_{i+1}\cdots a_n)\stackrel{(a_i,1)}\to 
(x, a_1\cdots a_n)$
of $\fF{K(X)}$
defining the commutative diagram in $K(X)$:
$$
\xymatrix{
x\cdot a_1\cdots a_n \ar[rr]^{1} & & x\cdot a_1\cdots a_n\\
x\cdot a_i \ar[u]^{a_1\cdots a_{i-1}a_{i+1}\cdots a_n} & & x \ar[ll]^{a_i} \ar[u]_{a_1\cdots a_n}
}
$$
The morphism $\Theta((\widehat{a_i,1}), \widetilde{x})$ is equal to morphism
$(x, a_1\cdots a_{i-1}a_{i+1}\cdots a_n)\stackrel{(1,a_i)}\to 
(x, a_1\cdots a_n)$ defining the commutative diagram in $K(X)$:
$$
\xymatrix{
x\cdot a_1\cdots a_{i-1}a_{i+1}\cdots a_n \ar[rr]^{a_i} & & x\cdot a_1\cdots a_n\\
x \ar[u]^{a_1\cdots a_{i-1}a_{i+1}\cdots a_n} & & x \ar[ll]^{1} \ar[u]_{a_1\cdots a_n}
}
$$
We get
\begin{gather*}
G_{\fF{K(S)}}\Theta^{op}((\widehat{1,a_i}), \widetilde{x})= 
G( (x\cdot a_i, a_1\cdots a_{i-1}a_{i+1}\cdots a_n)\stackrel{(a_i,1)}\to 
(x, a_1\cdots a_n))\\
G_{\fF{K(S)}}\Theta^{op}((\widehat{a_i,1}), \widetilde{x})= 
G( (x, a_1\cdots a_{i-1}a_{i+1}\cdots a_n)\stackrel{(1,a_i)}\to 
(x, a_1\cdots a_n))
\end{gather*}
Using Theorem \ref{main1}, we arrive to the formulas (\ref{maindiff}) for differentials. 
\hfill $\Box$

Let $\mC$ be an arbitrary small category. For any morphism 
$\alpha\in \mC(A,B)$, we denote by $\dom\alpha=A$ its domain and $\cod\alpha = B$ its codomain.
Extend the maps $\cod, \dom: \Ob\fF\mC \to \Ob\mC$ to morphisms by 
formulas $\cod(\alpha\stackrel{(f,g)}\to \beta)= g$ 
and $\dom(\alpha\stackrel{(f,g)}\to \beta)= f$. 
We have the functors  $\cod: \fF\mC\to \mC$ and $\dom: \fF\mC^{op}\to \mC$.

\begin{lemma}
For any object $c\in \Ob\mC$, the classifying space of the left fibre $(\cod/c)$ and right fibre $(c/\dom)$ 
homotopy equivalent to the single point. 
\end{lemma}
{\sc Proof.} By Quillen \cite{qui1973}, if there exists a natural transformation  $\eta: \Phi\to \Psi$
between functors $\Phi,\Psi: \mB\to \mD$, then the corresponding continuous maps 
of classifying spaces are homotopic.
Thus, if there are adjoint functors 
$\Phi: \mB\to \mD$  and $\Psi: \mD\to \mB$, then the classifying spaces 
of the categories $\mB$ and $\mD$ are homotopy equivalent. 
In \cite{X1997}, for a proof the sugestion \cite[Proposition 2.6]{X1997},
it was remarked that there exists a pair of adjoint functors
$
\xymatrix{(\mC/c)^{op}\ar@/_/[r] & \cod/c \ar@/_/[l]
}$.
Since $H_n(B\mC)\cong \coLim_n^{\mC}\Delta\ZZ$
\cite{qui1973}, this implies  $\coLim_n^{\cod/c}\Delta\ZZ=0$ for $n>0$, 
and  $\coLim_0^{\cod/c}\Delta\ZZ=\ZZ$.

Let us prove that the classifying space of the category $(c/\dom)$ homotopy equivalent 
to single point.
Any object in the category $(c/\dom)$ is a pair 
$(\alpha,\xi)$ of morphisms   
$c\stackrel{\xi}\to \dom\alpha \stackrel{\alpha}\to \cod\alpha$.
Morphisms $(\alpha, \xi)\to (\beta,\zeta)$ are given by commutative diagrams
$$
\xymatrix{
c \ar[rd]_{\zeta} \ar[r]^{\xi} & \dom\alpha\ar[d]_{f} \ar[r]^{\alpha} & \cod\alpha\\
& \dom{\beta} \ar[r]_{\beta} & \cod\beta\ar[u]_{g}
}
$$
Its full subcategory consisting of objects $(\alpha, 1_c)$, is isomorphic to
 $(c/\mC)^{op}$.
For every $(\zeta, \beta)\in \Ob(c/\dom)$, the following commutative diagram
$$
\xymatrix{
c \ar[rd]_{\zeta} \ar[r]^{1_c} & c\ar[d]_{\zeta} \ar[r]^{\beta\zeta} & \cod\beta\\
& \dom{\beta} \ar[r]_{\beta} & \cod\beta\ar[u]_{1}
}
$$
determines a morphism $(\beta\zeta, 1_c)\to (\beta,\zeta)$ such that for any another morphism
$(\alpha,1_c)\to (\beta,\zeta)$ there exists unique morphism $(\alpha, 1_c)\to (\beta,\zeta)$
making a commutative diagram
$$
\xymatrix{
(\beta\zeta, 1_c)\ar[rr] && (\beta,\zeta)\\
& (\alpha,1_c)\ar@{-->}[lu]^{\exists!} \ar[ru]
}
$$ 
This implies that the functor $U: (c/\mC)^{op}\to (c/\dom)$, $U(\alpha)=(\alpha,1_c)$ 
has the right adjoint $R(\beta,\zeta)=\beta\zeta$. 
Conclude that the classifying spaces of categories $ (c / \mC) $
and $(c/\dom)$ are homotopy equivalent, and hence $B (c/\dom)$ is homotopy equivalent to single point.
In particular, its $n$th homology groups equal $0$ for $n>0$, and $\coLim^{\dom/c}_0\Delta\ZZ= \ZZ$ for $n=0$.
\hfill $\Box$

We have proved that the integral homology groups of the fibers for functors 
$\fF\mC\stackrel{\cod}\to \mC$
and $\fF\mC\stackrel{\dom^{op}}\to \mC^{op}$ are isomorphic to the homology of single point.
Using Oberst's Theorem, as amended by \cite[Proposition 3.4]{X20081}, we obtain

\begin{proposition}
For any small category $\mC$, and the functor $F: \mC \to \Ab$ (respectively, $F:\mC^{op} \to \Ab$)
there exists an isomorphism $\coLim_n^{\mC} F \cong \coLim_n^{\fF\mC^{op}} (F \circ \dom)$
(respectively, $\coLim_n^{\mC^{op}}F \cong \coLim_n^{\fF\mC^{op}} (F \circ \cod^{op})$).
\end{proposition}

Let $L(S)$ denotes the free Abelian group generated by a set $S$.

\begin{corollary}\label{compcom}
Let $X$ be a set with right action of a locally finite trace monoid $M(E,I)$
and let  
$S\subseteq X$ be a subset for which  $K(S)\subseteq K(X)$ is the convex subcategory.
Then groups $\coLim_n^{K(S)}\Delta\ZZ$ will be isomorphic to the homology of the complex
\begin{multline*}
0 \leftarrow L(S) \stackrel{d_1}\leftarrow LQ_1(S,E,I) \stackrel{d_2}\leftarrow LQ_2(S,E,I)
\leftarrow \cdots \\
\cdots \leftarrow LQ_{n-1}(S,E,I) \stackrel{d_n}\leftarrow LQ_n(S,E,I) \leftarrow \cdots
\end{multline*}
with differentials 
\begin{multline*}
d_n(x, a_1, \cdots, a_n)= \sum\limits_{i=1}^n(-1)^i(x\cdot a_i, a_1, \cdots, a_{i-1}, a_{i+1}, \cdots, a_n)\\
			- \sum\limits_{i=1}^n(-1)^i(x, a_1, \cdots, a_{i-1}, a_{i+1}, \cdots, a_n)
\end{multline*}
\end{corollary}

\section{Applications}

Consider an algorithm for computing the integral homology groups
of a category associated with a finitely generated trace monoid action 
on a finite set. On this basis, we construct an algorithm for finding integral
homology groups of finite CE nets.

\subsection{Computing the homology groups of a category\\
 associated with a partial trace monoid action}

Let $M(E,I)$ be a trace monoid generated by a finite set  
$E$. Consider a 
 partial right action of $M(E,I)$ on a finite set $S$.
By Lemma \ref{KSclosed}, the category $K(S)$ is convex in $K_*(S)$. Hence, 
for computing the groups 
$\coLim_n^{K(S)}\Delta\ZZ$,  
we can use Corollary \ref{compcom}, substituting $X = US_*$.

Preliminary necessary define some arbitrary total order relation on $E$.

The algorithm consists of two steps.
 
First is a set of free Abelian groups
$C_n=\ZZ^{|Q_n(S,E,I)|}$ and matrix differentials
$d_1$, $d_2$, $\cdots$.
The differential $d_0=0$.
There is a correspondence between
elements $ (x, a_1, \cdots, a_n) \in Q_n (S, E, I) $ and the columns of the matrix of differential 
$d_n$. 
The strings correspond to the elements
in $Q_{n-1}(S,E,I)$. 
Since 
\begin{multline*}
d_n(x, a_1, \cdots, a_n)= \sum\limits_{i=1}^n(-1)^i(x\cdot a_i, a_1, \cdots, a_{i-1}, a_{i+1}, \cdots, a_n)\\
			- \sum\limits_{i=1}^n(-1)^i(x, a_1, \cdots, a_{i-1}, a_{i+1}, \cdots, a_n),
\end{multline*}
we fill the column corresponding to the element $(x, a_1, \cdots, a_n)$ as follows:
for every $i\in\{1, \cdots, n\}$, we would write $(-1)^i$ in the string
 corresponding to 
$$
(x\cdot a_i, a_1, \cdots, a_{i-1}, a_{i+1}, \cdots, a_n)
$$
 and write $(-1)^{i+1}$ in the string
$$
(x, a_1, \cdots, a_{i-1}, a_{i+1}, \cdots, a_n).
$$ 
Other entries of the matrix equal $0$.
For each $ n \geq 1 $, matrix of  $ d_n $ to be transformed into
the Smith normal form (see \cite {X20042}), consisting of the diagonal matrix
with entries $ (\delta^{n}_1, \delta^{n}_2, \cdots, \delta^{n}_{rank (d_{n})})$.
The homology groups are equal
$$
H_n= 
\ZZ^{|Q_n(S,E,I)|-rank(d_n)-rank(d_{n+1})}\oplus 
\ZZ/\delta^{n+1}_1\ZZ \oplus \cdots \oplus \ZZ/\delta^{n+1}_{rank(d_{n+1})}\ZZ,
$$
where $\ZZ / \delta^{n +1}_k \ZZ $ denote the group of residues modulo $ \delta^{n +1}_k $

For example, we consider the trace monoid $M(E,I)$ with 
$E=\{a_1, a_2, a_3\}$ and $I=\{(a_i, a_j)| 1\leq i, j \leq 3 ~\&~ i\not=j\}$.
In this case, $M(E,I)$ will be the free commutative monoid generated 
by three elements. Define the action on the set
$S= \{s_0, s_1, s_2, s_3, s_4, s_5, s_6, s_7\}$
specifying the elements
 $s_i\cdot a_j$, $0\leq i\leq 7$, $1\leq j\leq 3$
by the graph with labels consisting of the edges
 $s_i\stackrel{a_i}\to s_i\cdot a_j$:

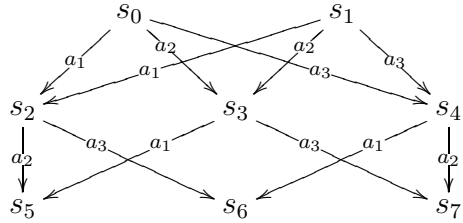
$$
\xymatrix{
 &	s_0\ar[ld]|{a_1} \ar[rd]|(.35){a_2} \ar[rrrd]|(.6){a_3}	
& & s_1	
 \ar[llld]|(.6){a_1} \ar[ld]|(.35){a_2}  \ar[rd]|{a_3}	
\\
s_2 \ar[d]|{a_2} \ar[rrd]|(.35){a_3}	&  
&	s_3 \ar[lld]|(.35){a_1} \ar[rrd]|(.35){a_3}	& &	
s_4 \ar[lld]|(.35){a_1} \ar[d]|(.5){a_2}\\
s_5	& &	
s_6	& &	s_7
}
$$
The set $Q_0(S,E,I)=S$ consists of 8 elements,
\begin{multline}\label{enum1}
Q_1(S,E,I)=\\
 \{(s_0,a_1), (s_0,a_2), (s_0,a_3), (s_1,a_1), (s_1,a_2), (s_1,a_3), \\
(s_2,a_2), (s_2,a_3), (s_3,a_1), (s_3,a_3), (s_4,a_1), (s_4,a_2)\}
\end{multline}
\begin{multline*}
Q_2(S,E,I)=\\
 \{(s_0,a_1, a_2), (s_0,a_1, a_3), (s_0, a_2, a_3), 
(s_1,a_1, a_2), (s_1,a_1, a_3), (s_1,a_2, a_3)\}
\end{multline*}
We have $|S|=8$, $|Q_1(S,E,I)|=12$, $|Q_2(S,E,I)|=6$, and $|Q_n(S,E,I)|=0$ for $n>2$.
Hence, the complex consists of Abelian groups
$$
	0 \stackrel{d_0}\leftarrow \ZZ^8 \stackrel{d_1}\leftarrow 
\ZZ^{12} \stackrel{d_2}\leftarrow \ZZ^6 
\stackrel{d_3}\leftarrow 0 \leftarrow \cdots 
$$
The entries of the matrix for $d_1$ are found from the formula
\begin{equation}\label{fdiff1}
d_1(s,a)=-s\cdot a+s.
\end{equation} 
Strings of the matrix will correspond to the elements of $S$.
The first line to $ s_0 $, second to $ s_1 $,  and so on.
 The first column corresponds to $ (s_0, a_1)$,
second to $ (s_0, a_2) $ and so on, according to the formula (\ref{enum1}).
The entries of the first column we find by the formula 
$d_1 (s_0, a_1) =-s_0 \cdot {a_1} + s_0 =-s_2 + s_0$.
Hence, in the string corresponding to the $ s_2 $, we must write  $-1$, and in line $s_0$, we must write $+1$.
For the second column $ d_1 (s_0, a_2) =-s_0 a_2 + s_0 =-s_3 + s_0 $. Write down in the fourth string
the second column of $ -1 $, and the first - $ 1 $. Go to the next column, etc.
We will obtain the matrix
$$
\left(
\begin{array}{rrrrrrrrrrrr}
1 & 1 & 1 & 0 & 0 & 0 & 0 & 0 & 0 & 0 & 0 & 0 \\
0 & 0 & 0 & 1 & 1 & 1 & 0 & 0 & 0 & 0 & 0 & 0 \\
-1& 0 & 0 &-1 & 0 & 0 & 1 & 1 & 0 & 0 & 0 & 0 \\ 
0 &-1 & 0 & 0 &-1 & 0 & 0 & 0 & 1 & 1 & 0 & 0 \\
0 &0  &-1 & 0 & 0 &-1 & 0 & 0 & 0 & 0 & 1 & 1 \\
0 & 0 & 0 & 0 & 0 & 0 &-1 & 0 &-1 & 0 & 0 & 0 \\
0 & 0 & 0 & 0 & 0 & 0 & 0 &-1 & 0 & 0 & -1& 0 \\
0 & 0 & 0 & 0 & 0 & 0 & 0 & 0 & 0 &-1 & 0 & -1  
\end{array}
\right)
$$
The matrix would be transformed
to the Smith diagonal form whose all entries are equal to zero,
except 
$a_{11}= a_{22}= a_{33}= a_{44}= a_{55}= a_{66}= a_{77}= 1$.
This implies
$$
\coLim^{K(S)}_0{\Delta\ZZ}= 
\ZZ^{8-rank(d_0)-rank(d_1)}
\oplus \underbrace{\ZZ/1\ZZ\oplus\cdots \oplus\ZZ/1\ZZ}_{7~\mbox{times}}
=\ZZ^{8-0-7}=\ZZ.
$$

Write entries of  $d_2$ by formula
\begin{equation}\label{fdiff2}
d_2(s,e_1,e_2)= -(s\cdot e_1, e_2)+(s\cdot e_2, e_1) + (s, e_2) - (s, e_1).
\end{equation}
First column corresponds to $(s_0, a_1, a_2)$. We have 
\begin{multline*}
d_2(s_0, a_1, a_2)= \\
-(s_0\cdot a_1, a_2)+ (s_0\cdot a_2, a_1)+ (s_0, a_2)- (s_0, a_1)
= \\
-(s_2, a_2)+ (s_3, a_1)+ (s_0, a_2) - (s_0, a_1).
\end{multline*}
In the strings corresponded to $(s_2, a_2)$ and $(s_0, a_1)$, we write $-1$, 
but in the strings $(s_3,a_1)$ and $(s_0,a_2)$ we write $+1$.
Similarly, the entries are written for other columns. The matrix of $d_2$
equals
$$
\left(
\begin{array}{rrrrrr}
-1 & -1 & 0 & 0 & 0 & 0  \\
1 & 0 & -1 & 0 & 0 & 0  \\
0 & 1 & 1 & 0 & 0 & 0  \\
0 & 0 & 0 & -1 & -1 & 0  \\
0 & 0 & 0 & 1 & 0 & -1  \\
0 & 0 & 0 & 0 & 1 & 1  \\
-1 & 0 & 0 & -1 & 0 & 0  \\
0 & -1 & 0 & 0 & -1 & 0  \\  
1 & 0 & 0 & 1 & 0 & 0  \\
0 & 0 & -1 & 0 & 0 & -1  \\
0 & 1 & 0 & 0 & 1 & 0  \\
0 & 0 & 1 & 0 & 0 & 1    
\end{array}
\right)
$$
Its normal form equals to the matrix with $a_{11}= a_{22}= a_{33}= a_{44}= a_{55}=1$ 
and other entries are $0$.
We obtain
$$
\coLim_1^{K(S)}\Delta\ZZ= \ZZ^{12-rank(d_1)-rank(d_2)}
\oplus \underbrace{\ZZ/1\ZZ \oplus \cdots \oplus \ZZ/1\ZZ}_{5~\mbox{times}} =\ZZ^{12-7-5}=0
$$
Finally, $d_3$ is null matrix. Find 
$$
\coLim_2^{K(S)}\Delta\ZZ= \ZZ^{6-rank(d_2)-0}= \ZZ.
$$

\subsection{Computing the homology groups of CE nets}

We shall use the terminology of paper \cite{X20042} in which 
homology groups of a CE net have been introduced by homology of corresponding 
asynchronous system.

By  \cite{X20042}, an {\em asynchronous system} can be defined as a triple
$(S, s_0, M(E,I))$ consisting of a partial trace monoid action 
of $M(E,I)$ on a set $S$ with a distinguished element $s_0\in S$.
Elements $s\in S$ are called {\em states} and elements of
$S(s_0)= \{s\cdot\mu| \mu\in M(E,I)\} \subseteq S$ are {\em reachable states}.
{\em Homology groups of asynchronous system with coefficients in 
a functor  $F: K(S)\to \Ab$} 
are Abelian groups
$\coLim_n^{K(S(s_0))}F|_{K(S(s_0))}$.

For a set $B$, denote by $2^B$ the set of all its subsets.

A {\em CE net} \cite{X20042} or {\em Petri net} \cite{win1995} 
 is a quintuple $(B,E, pre, post, s_0)$ consisting of finite sets
$B$ and $E$, the maps $pre, post: E\to 2^{B}$  
and a subset $s_0\subseteq B$. Here CE is abbreviation for the words Conditions/Events.

Let ${\mathcal N} =(B,E, pre, post, s_0)$ be a CE net.
Define an relation $I\subseteq E\times E$ as the set of pairs $(a,b)$ 
for which $(pre(a)\cup post(a))\cap(pre(b)\cup post(b))=\emptyset$.
To every element $e \in E $ we assign a partial mapping $(-)\cdot{e}: 2^B \rightharpoonup 2^B$
defined for $s \subseteq B$ satisfying to the condition
$$
(pre(e)\subseteq s) \quad \& \quad (post(e)\cap s=\emptyset).
$$ 
In these cases, we take $s\cdot{e}=(s\setminus pre(e)) \cup post(e)$ \cite{maz1987}.
This define a partial action of $M(E,I)$ on the set $2^{B}$. 
Assuming $S=2^B$, we get an asynchronous system $(S, s_0, M(E,I))$, which corresponds to 
the CE net ${\mathcal N}= (B,E, pre, post, s_0)$.
The homology groups  $H_n({\mathcal N})$ a defined as 
$\coLim_n^{K(S(s_0))}\Delta\ZZ$ where $S(s_0)$ is the set of reachable states.

The problem of computing the homology groups of CE nets
is solved by the method described above, applied to a partial trace monoid action 
of $M(E,I)$ on set $ S (s_0) $ by the formula
$s\mapsto s\cdot{e}$ for $s\in S$ and $e\in E$.

For example, consider the following  CE net (of a {\em pipeline})

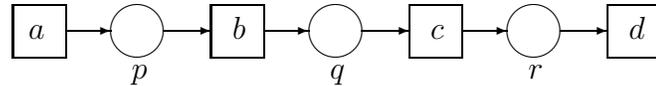
\begin{center}
\begin{picture}(320,40)
\multiput(97,20)(75,0){3}
{\circle{20}}

\put(95,2){$p$}
\put(170,2){$q$}
\put(245,2){$r$}

\multiput(50,10)(75,0){4}
{\line(1,0){20}}
\multiput(70,10)(75,0){4}
{\line(0,1){20}}
\multiput(70,30)(75,0){4}
{\line(-1,0){20}}
\multiput(50,30)(75,0){4}
{\line(0,-1){20}}

\multiput(70,20)(75,0){3}
{\vector(1,0){17}}
\multiput(107,20)(75,0){3}
{\vector(1,0){18}}

\put(56,17){$a$}
\put(133,17){$b$}
\put(208,17){$c$}
\put(283,17){$d$}

\end{picture}
\end{center}
The corresponding asynchronous system $(S, s_0, M(E,I))$ will consist 
of the set $S=2^{\{p,q,r\}}$, 
$s_0=\emptyset$, $E=\{a,b,c,d\}$, and $I=\{(a,c), (a,d), (b,d)\}$. 
Subsets $s\subseteq\{p, q, r\}$ will be given as triples
$(\varepsilon_p, \varepsilon_q, \varepsilon_r)$, 
where $\varepsilon_x=1 \Leftrightarrow x\in s$.
For example, the triple $(1,0,1)$ denotes the state $s=\{p,r\}$.

Actions of elements $e\in E$ are illustrated by the labeled directed graph 
(\ref{pipelineact})  with the arrows $s\stackrel{e} \to s'$
specified in the case $s \cdot e = s'$.

\begin{equation}\label{pipelineact}
\xymatrix{
(0, 1, 1)\ar[ddd]_a\ar[rrr]^d &&& (0,1,0)\ar[lld]_c \ar[ddd]^{a}\\
& (0,0,1) \ar[d]^a \ar[r]_d & (0,0,0)\ar[d]_a\\
& (1,0,1)\ar[luu]^b \ar[r]^d & (1,0,0)\ar[ruu]_b\\
(1, 1, 1)\ar[rrr]_d &&& (1,1,0)\ar[llu]^c
}
\end{equation}
We have $S(s_0)=S$. The set $S(s_0)$ has $8$ elements,
\begin{multline*}
Q_1(S,E,I)= \{((0,0,0),a), ((0,0,1),a), ((0,0,1),d), ((0,1,0),a),\\
((0,1,0),c),((0,1,1),a),((0,1,1),d), ((1,0,0),b),\\
((1,0,1),b),
((1,0,1),d),((1,1,0),c),((1,1,1),d)
\}
\end{multline*}
$$
Q_2(S,E,I)= \{((0,0,1),a,d), ((0,1,0),a,c), ((0,1,1),a,d), ((1,0,1),b,d) \}
$$
The complex for the computing the groups $H_n({\cal N})$ consists of the Abelian groups
$$
0 \stackrel{d_0}\leftarrow \ZZ^{8} \stackrel{d_1}\leftarrow \ZZ^{12} \stackrel{d_2}\leftarrow \ZZ^{4}
		\stackrel{d_3}\leftarrow 0 \leftarrow \cdots .
$$
Entries of $d_1$ compute by the formula  (\ref{fdiff1}).
$$
\left(
\begin{array}{rrrrrrrrrrrr}
1 & 0 & -1 & 0 & 0 & 0 & 0 & 0 & 0 &0 & 0 & 0 \\
0 & 1 &  1 & 0 &-1 & 0 & 0 & 0 & 0 &0 & 0 & 0 \\
0 & 0 &  0 & 1 & 1 & 0 &-1 &-1 & 0 &0 & 0 & 0 \\
0 & 0 &  0 & 0 & 0 & 1 & 1 & 0 &-1 &0 & 0 & 0 \\
-1 & 0 & 0 & 0 & 0 & 0 & 0 & 1 & 0&-1 & 0 & 0 \\
0 & -1 & 0 & 0 & 0 & 0 & 0 & 0 & 1 &1 &-1 & 0 \\
0 & 0 &  0 &-1 & 0 & 0 & 0 & 0 & 0 &0 & 1 &-1 \\
0 & 0 &  0 & 0 & 0 &-1 & 0 & 0 & 0 &0 & 0 & 1 
\end{array}
\right)
$$
The Smith normal form of $d_1$ will be consists of $7$ unities.
Entries of the matrix $d_2$ find by the formula (\ref{fdiff2}):
$$
\left(
\begin{array}{rrrr}
 1 & 0 & 0 & 0\\
-1 & 1 & 0 & 0\\
 1 & 0 & 0 & 0\\
 0 & -1 & 1 & 0\\
 0 & 1 & 0 & 0\\
 0 & 0 & -1 & 0\\
 0 & 0 & 1 & -1\\
 0 & 0 & 0 & 1\\
 0 & 0 & 0 & -1\\
-1 & 0 & 0 & 1\\
 0 &-1 & 0 & 0\\
 0 & 0 & -1 & 0
\end{array}
\right)
$$
The Smith normal form consists of $4$ unities.
We obtain
$H_0({\mathcal N})=\ZZ^{8-0-7}=\ZZ$, 
$H_1({\mathcal N})=\ZZ^{12-7-4}=\ZZ$,
$H_2({\mathcal N})=\ZZ^{4-4-0}=0$.
Consequently, $0$th and $1$th homology groups 
of the CE net are equal to $\ZZ$. 
Other homology groups are $0$.

\medskip
{\bf Acknowledgements}. I would like to express my gratitude for financial support of this research by
Komsomolsk-on-Amur State Technical University. The research was also supported by T\"UBITAK and
NATO in Turkey and Grant Center at Novosibirsk State University.

\end{document}